\documentclass[10pt]{article}

\usepackage{graphics}
\usepackage{graphicx}

\usepackage[a4paper, left=35mm,right=35mm,top=34mm,bottom=34mm]{geometry}
\usepackage[utf8]{inputenc}
\usepackage[T1]{fontenc}
\usepackage[english]{babel}

   \graphicspath{{img/pdf/}{img/jpeg/}{img/png/}}
   \DeclareGraphicsExtensions{.pdf,.jpg,.jpeg,.png}
\usepackage{0_packages/akca}
\usepackage{xspace}

\author{
  {\normalsize Abal-Kassim Cheik Ahamed}\thanks{CUDA Research Center, Applied Mathematics and Systems Laboratory, CentraleSup\'elec, Universit\'e Paris-Saclay, France.}
	\and
  {\normalsize Fr\'ed\'eric Magoul\`es}\thanks{CUDA Research Center, Applied Mathematics and Systems Laboratory, CentraleSup\'elec, Universit\'e Paris-Saclays, France
    (correspondence, frederic.magoules@hotmail.com).}
		}	
\title{On the stability and performance of the solution of sparse linear systems by partitioned procedures}
\date{}

\begin{document}
\maketitle
\thispagestyle{fancy}

\begin{abstract}
\noindent In this paper, we present, evaluate and analyse the performance of parallel synchronous Jacobi algorithms by different partitioned procedures including band-row splitting, band-row sparsity pattern splitting and substructuring splitting, when solving sparse large linear systems.
Numerical experiments performed on a set of academic 3D Laplace equation and on a real gravity matrices arising from the Chicxulub crater are exhibited, and show the impact of splitting on parallel synchronous iterations when solving sparse large linear systems.
The numerical results clearly show the interest of substructuring methods compared to band-row splitting strategies.
\end{abstract}


\section{Introduction}
\label{sec:introduction}

In applied sciences a system of partial differential equations (PDEs) is used as a fundamental physical and mathematical model which is solved numerically. 
To obtain a numerical solution of PDE, discretization methods (finite difference/volume/element methods) are applied, and linear systems with very large sparse matrices are obtained after some linearization process, if necessary. 
In scientific computing, solving large linear systems is often the most expensive part, both in terms of memory space and computation time.
A system of linear equations can be written in the matrix notation as
\begin{equation}
\label{eq:linear_equation}
A x = b,
\end{equation}
where $A$ is a $n \times n$ matrix, $b$ the right-hand side in $\C^n$ and $x$ in $\C^n$ represents the solution vector we are looking for.
The system~\eqref{eq:linear_equation} has a unique solution if and only if the matrix $A$ is non-singular, which means that the determinant of the matrix $A$ is non-zero.
There are several methodologies to solve this type of large linear system, \eg, sparse direct solvers, iterative solvers, and a combination of those.
In this paper, we focus on parallel synchronous Jacobi methods with different splitting strategies to efficiently solve large and sparse linear systems, arising from the discretization of finite element methods.
The synchronous parallel algorithms have been for solving these systems in a context of heterogeneous multi-core systems.\\

The rest of the paper is structured as follows.
Section~\ref{sec:Background in Parallel Computing} introduces the motivations and backgrounds of parallel computing in scientific computing.
In Section~\ref{sec:Jacobi iterations}, we describe Jacobi algorithm and present the main points of its implementation in parallel for different partitioning, leading to the design of linear algebra operations.
We also give some properties and convergence results of Jacobi splitting methods.
Section~\ref{sec:Matrix partitioning} presents the partitioning procedures of the data.
This section also details substructuring methods.
The next section (Section~\ref{sec:Benchmark environments}) briefly presents the kernel of the library we used for the computation of linear algebra operations. The hardware configuration used in the experiments is also described in Section~\ref{sec:Benchmark environments}.
Numerical experiments and test cases are reported and discussed in Section~\ref{sec:Numerical results}.
Finally, Section~\ref{sec:conclusion} concludes this paper.

\section{Background in Parallel Computing}
\label{sec:Background in Parallel Computing}

  \subsection{Parallel architectures}

Parallel computing is a set of hardware and software techniques allowing the simultaneous execution of sequences of independent instructions on different processes and/or cores.
The set of hardware and software techniques consists of various architectures of parallel computers and various models of parallel programming.

Parallel computing~\cite{batcher_designmassivelyparallel_1980}~\cite{leighton_introductionparallelalgorithms_1992}~\cite{foster_designingbuildingparallel_1995}~\cite{ciegis_parallelscientificcomputing_2008} has two objectives: accelerate the execution of a program by distributing work and executing a big problem with more material resources (memory, storage capacity, processor, etc.).
A parallel computer can be: a multi-core processor having at least two physical computing units on the same chip or a supercomputer, which gathers the components of several computers (processors and memories) in only one machine, or a distributed platform made up of several independent machines, homogeneous or heterogeneous, connected with one another by a communication network.
In the jargon, the supercomputer is more often called a cluster.\\

  \subsubsection{Classifications of parallel architectures}

In literature, there are several classifications of parallel architectures based on various criteria.
In this paper, we present briefly the most largely used classification in parallel computing: the \emph{taxonomy of Flynn}~\cite{flynn_computerorganizationstheir_1972}.
This classifies parallel architectures in several categories, according to whether great volumes of data are processed in parallel (or not) or whether a large number of instructions are being carried out at the same time (or not).
In the taxonomy of Flynn, one distinguishes four classes: \emph{SISD (Single Instruction, Single Data)}, \emph{SIMD (Single Instruction, Multiple Data)}, \emph{MISD (Multiple Instruction, Single Data)} and \emph{MIMD (Multiple Instruction, Multiple Data)}.
The \emph{SISD} class represents most conventional computing equipment which has only one computing unit.
The sequential computer (traditional uniprocessor machine) can process only one operation at a time.
The \emph{SIMD} class corresponds to computers which have a large number of computing units, \eg array processors or GPUs.
The computer can exploit multiple data against a single instruction to perform operations which may be naturally parallelized.
At each clock cycle, all the processors of a SIMD computer carry out the same instruction on different data simultaneously.
The \emph{MISD} class corresponds to the parallel machines which perform several instructions, simultaneously, on the same data.
These architectures are uncommon and are generally used for fault tolerance.
\emph{MIMD} represents the most used class in the taxonomy of Flynn.
In this architecture, multiple autonomous processors simultaneously execute different instructions on different data, \eg distributed systems.
They exploit a single shared memory space or a distributed memory space.
Computers using \emph{MIMD} have a number of processors that function asynchronously and independently.

  \subsubsection{Memory model}

In parallel computing, we have two main models of memory management~\cite{carlini_transputersparallelarchitectures_1990}~\cite{leighton_introductionparallelalgorithms_1992}: shared memory and distributed memory.
These two models help define the mode of data access to the data of the other cooperating processors in a parallel computation for a given process.
In a shared memory model, the processors are all connected to a "globally available" memory.
All the processors have a direct access to the same physical memory space via powerful links of communication (\eg cores of a single machine).
The second memory model is generally used in parallel computers with distributed resources, \eg a grid or cluster.
In this model, each processor has its own individual memory location.
Each processor has no direct knowledge of other processors' memory.
However, access to the data of the distributed memories is ensured by message passing between cooperating processors through a communication network.
In this paper, we focus on algorithms in parallel distributed-memory computing.

  \subsection{Parallel algorithms}

Parallel computations are fundamental and ubiquitous in numerical analysis and large application areas, when we deal with the problem of large-size.
Distributed computing constantly gains in importance and has become an important tool in common scientific research work~\cite{bertsekas_paralleldistributedcomputation_1989}.
When computation is done in parallel, the simplest solution consists in synchronizing the processors during iterations.
In that particular field of parallel programming, the commonly used model is the synchronous message passing.
In the synchronous case~\cite{bahi_performancecomparisonparallel_2004}~\cite{bahi_efficientrobustdecentralized_2008}, communications are strongly penalizing the overall performance of the application.
Indeed, they often involve large idle periods, especially when the processors are heterogeneous.
The performance of algorithms strongly depends on the management of interprocessor communication.
Parallel algorithms require that a large number of parameters, such as the number and the typology of the processors, be taken into account.
Therefore, in parallel processing, the step of data distribution is crucial and strongly impacts the performance of algorithms.\\

In this paper, we pay special attention to substructuring method, which is the precursor of non-overlapping domain decomposition methods~\cite{letallec_domaindecompositionmethods_1994}~\cite{smith_domaindecompositionparallel_1996}~\cite{quarteroni_domaindecompositionmethods_1999}~\cite{toselli_domaindecompositionmethods_2005}\cite{gander_analysispatchsubstructuring_2007}~\cite{magoules_substructuringtechniquesdomain_2010}.

\section{Jacobi iterations}
\label{sec:Jacobi iterations}

The Jacobi method~\cite{margaris_parallelimplementationsjacobi_2014} is a stationary iterative method based on splittings~\cite{oleary_multisplittingsmatrices_1985} of the matrix $A \in \K^{n \times n}$.
Let
\begin{equation}
\label{eq:substructuring_splitting}
A = M - N,
\end{equation}
be a splitting of $A$, such that $M$ is a non-singular matrix.
The iterative algorithm associated with the splitting~\eqref{eq:substructuring_splitting} is defined by
\begin{equation}
\label{eq:iterative_system_algorithm}
  \begin{cases}
    \uk{0} \text{ given,} \\
    \uk{k+1} = M^{-1}N\uk{k} + M^{-1}b, \quad \forall\ k \in\N.
  \end{cases}
\end{equation}
Let $T = M^{-1}N$ and $c = M^{-1}b$. The matrix $\abs{T}$ is defined such that $\abs{T} = (\abs{T_{i,j}})$, for $i,j \in \intervaldot{1}{n}$.
The speed of the algorithm is determined by the quantity $\rho(T)$, where $\rho(T)$ denotes the spectral radius of the matrix $T$.
More accurately, the bounded value of the error between the exact solution and the approximate solution at iteration $k$ is expressed by:
\begin{equation}
  \label{eq:error_lim_exact_solution_and_approximate_solution}
  \forall \mathepsilon>0, \quad \norm{\uk{k}-u^*} \leq (\rho(T)+\mathepsilon)^k \norm{\uk{0}-u^*}
\end{equation}
\begin{theorem}\thmlabel{jacobi_conv_seq_sync}
The sequential and parallel synchronous algorithms~\eqref{eq:iterative_system_algorithm} converge if and only if the spectral radius $\rho(T) < 1$.
\end{theorem}
The Jacobi algorithm is obtained by considering $M$ and $N$ in algorithm~\eqref{eq:iterative_system_algorithm} such as $M = D$ is the diagonal part of $A$ (non-singular), $N = -(L + U)$ where $L$ and $U$ are the strictly lower and upper triangular of $A$, \ie
\begin{equation}
\label{eq:jacobi_system_iteration}
  \begin{cases}
    \uk{0} \text{ given,} \\
    \uk{k+1} = D^{-1} \left( b - (L+U)\uk{k} \right), \quad \forall\ k \in\N
  \end{cases}
\end{equation}
  \begin{equation}
  \label{eq:rewritten_jacobi_stationary_iterative_method}
  \begin{cases}
    \uk{0}_{i} \text{ given,} \\
    \uk{k+1}_{i} = \dfrac{1}{a_{ii}} \left( b_{i} - \sum_{\substack{j=1\\j \neq i}}^{n} a_{ij}\uk{k}_{j} \right),\;\;\; i = \intervaldot{1}{n},\;\;k \in \N.
  \end{cases}
  \end{equation}
\figreft{design_jacobi_algorithm} shows the scheme of the Jacobi method~\eqref{eq:jacobi_system_iteration}.
\includefig{
\centering
\scalebox{0.65}{
\LARGE
\begin{tabular}{m{5cm}m{1.6cm}m{0.6cm}m{5cm}m{1.5cm}m{0.5cm}m{0.5cm}}
  \begin{tikzpicture}
  \draw[color=black,line width=0.02cm,fill=white,fill opacity=1] (0,0) rectangle (5,5);
  \draw[color=black,line width=0.02cm,fill=gray,fill opacity=1] (4.5,0) -- (5,0) -- (5,0.5) -- (0.5,5) -- (0,5) -- (0,4.5) -- (4.5,0);
  \draw[color=black] (2.5,2.5) node {$D$};
  \end{tikzpicture}
& $\uk{k+1}$ & + &
  \begin{tikzpicture}
  \draw[color=black,line width=0.02cm,fill=white,fill opacity=1] (0,0) rectangle (5,5);
  \draw[color=black,line width=0.02cm,fill=gray,fill opacity=1] (0,0) -- (0,4.5) -- (4.5,0) -- (0,0);
  \draw[color=black,line width=0.02cm,fill=gray,fill opacity=1] (5,0.5) -- (0.5,5) -- (5,5) -- (5,0.5);
  \draw[color=black] (3.8,4) node {$U$};
  \draw[color=black] (1.2,1) node {$L$};
  \end{tikzpicture}
& $\uk{k}$ & = & $b$\\
\end{tabular}
} 
}{Design of the Jacobi algorithm}{}{design_jacobi_algorithm}
Considering the algorithm~\eqref{eq:rewritten_jacobi_stationary_iterative_method}, the computations of the components $\uk{k}_{i}, i \in \intervaldot{1}{n}$ are independent, which signifies that their update can be performed in parallel.
The Jacobi algorithm converges more slowly, but is \textit{fully parallelizable}.
The jacobi iteration given in algorithm~\eqref{eq:jacobi_system_iteration} can be rewritten in the following form,
\begin{equation}
\label{eq:jacobi_system_iteration_2}
  \begin{cases}
    \uk{0} \text{ given,} \\
    \uk{k+1} = D^{-1}\left( b - A\uk{k} \right) + \uk{k}, \quad \forall\ k \in\N
  \end{cases}
\end{equation}
as described by the authors in~\cite{cheikahamed_efficientimplementationjacobi_2016}.
\algorithmreft{jacobi-general_algorithm_vectorial} present \emph{vectorial version} of Jacobi algorithm~\cite{cheikahamed_efficientimplementationjacobi_2016}.
The parallel version consists in computing all operations locally and then exchange the local residual between cooperating processors in order to compute the global convergence.
The crucial points in parallel iterative algorithms consists of the stoping criteria and the convergence detection steps.
The stopping criterion implemented in this work is based on \emph{simultaneous local convergence}~\cite{bahi_decentralizedconvergencedetection_2005}~\cite{bahi_synchronousasynchronoussolution_2006}~\cite{bahi_efficientrobustdecentralized_2008} of all the processors.
The convergence of the Jacobi method is ensured if the matrix $A$ is diagonally dominant: $\forall\ i \in \intervaldot{1}{\p}, \abs{a_{ii}} \ge \sum_{\substack{j=1\\j \neq i}} \abs{a_{ij}}$.
The authors presented in~\cite{cheikahamed_efficientimplementationjacobi_2016}, an original Jacobi implementation that helps optimize the solution of sparse linear systems on GPU-based implementation.
\includealgo[
  \Input{$\;n$: size of the matrix, \\
         $\;A$: $n \times n$ square matrix,
         $\;b$: right-hand side vector,\\
         $\;u^{(0)}$: initial guess,\\
         $\;\mathepsilon$: tolerance threshold,
         $\;K$: maximum number of iterations}
  \Output{$\;u$: solution vector}
]{
  Choose an initial guess $\uk{0}$ to the solution\\
  $k \gets 0$\\
 \emph{Compute} $D^{-1}$ \emph{//-~- $D^{-1}$ is computed once before the iterations}\algonllabel{jacobi_vectorial_inverseD}\\
  \While{Loop until convergence} {
    $q \gets A \uk{k}$\\
    \tcp*[l]{Compute $r \gets b - A \uk{k}$}
    $r \gets b - q$\\
   \emph{//-~- Compute $\norm{r} \equiv \norm{\uk{k+1}-\uk{k}}/\norm{D^{-1}}$}\\
    $\norm{r} \gets \norm{b - A \uk{k}}$\algonllabel{jacobi_vectorial_norm}\\
    \If{$\norm{r} \le \mathepsilon$} {
      $k \gets k + 1$\\
      break;
    }
    $\uk{k+1} \gets \uk{k} + D^{-1} r$\\
    $k \gets k + 1$
  }
}{Jacobi method: vectorial version}{}{jacobi-general_algorithm_vectorial}

\section{Matrix partitioning}
\label{sec:Matrix partitioning}

The main step in parallel processing consists in distributing the data on the cluster processors, which is commonly called parallel distributed computing.
In this section, we describe how data are distributed among processors for different splitting strategies: band-row, band-column, and substructuring splitting.
The distribution of data is accomplished as a preprocessing step, independently from the solver code.
The data such as matrix, right hand-size, vector solution and local to global, are written into a file, and will be input for the solver code.
The matrix is read from and written into the matrix market file~\cite{boisvert_matrixmarketweb_1997}~\cite{georgiev_numericalexperimentsapplying_2012}.\\

In this work, we have developed and implemented a code for partitioning data in order to control and adapt the splitting to the studied parallel algorithms.
This choice has allowed us to be free and implement parallel algorithms as we like.
The substructuring splitting uses METIS software~\cite{karypis_metissoftwarepackage_1998} for partitioning graphs and then perform matrix partitioning into sub-structures.\\

To illustrate the advantages and disadvantages of these types of partitioning, we will focus on the matrix-vector product, which is the most expensive operation in linear algebra~\cite{cheikahamed_fastsparsematrix_2012}~\cite{magoules_alineaadvancedlinear_2015}.

  \subsection{Sparse matrix formats}
The treatment of large sparse matrices in parallel~\cite{hassani_analysissparsematrix_2013} requires a good choice of storage format that helps the computations of involved operations.
The basic idea behind sparse matrix storage is to store only the non-zero matrix elements.
The distribution of non-zero coefficients depends on the features of the original problem.
The performance of the algorithms strongly depends on the data structure of the sparse matrices as demonstrated in~\cite{bolz_sparsematrixsolvers_2003}~\cite{xu_sparsematrixvector_2010}~\cite{bahi_parallelgmresimplementation_2011}~\cite{ren_sparselufactorization_2012}~\cite{suchoski_adaptingsparsetriangular_2012}~\cite{magoules_autotunedkrylov_2014}~\cite{magoules_fastiterativesolvers_2015}.\\

In this work, the matrices are stored in compressed sparse row (CSR) format in order to
optimize the memory storage and to make advantage of sparse structure for memory access.
Compressed Sparse Row (\emph{CSR}), described in~\figreft{csr_format}, is widely used because of minimal memory usage and the simplicity of the implementation.

To illustrate this storage formart, let us consider the sparse matrix $A$ described in \figsubreft{matrix_example}. In the following example, we consider a matrix indexed from $1$.
\figsubreft{matrix_pattern} draws the pattern of non-zero values of the matrix $A$.
\includesubfig{
\raisebox{10mm}{
  \scalebox{1.0}{
    $A=\begin{pmatrix}
    \colorbox{indiagreen}{-5} & \colorbox{lightgray}{14} & 0 & 0 & 0 \\
     0 & \colorbox{indiagreen}{8} & \colorbox{lightgray}{1} & 0 & 0 \\
     \colorbox{indiagreen}{2} & 0 & \colorbox{lightgray}{10} & 0 & 0 \\
     0 & \colorbox{indiagreen}{4} & 0 & \colorbox{lightgray}{2} & \colorbox{lightgray}{9} \\
     0 & 0 & \colorbox{indiagreen}{15} & 0 & \colorbox{lightgray}{7}
    \end{pmatrix}$
  }
}
}{Example matrix $A$, \colorbox{indiagreen}{nnz} first non-zero on the row}{matrix_example}{
  \scalebox{0.6}{
    \SetDisplayColumnIndexD{1}{0}
    \SetDisplayRowIndexD{1}{0}
    \SetColorFirstNnz{indiagreen}
    \SetSymbolDefaultNnz{$\bigstar$}
    \SetSymbolFirstNnz{$\blacksquare$}
    \def\nnzcoef{1/1,1/2,2/2,2/3,3/1,3/3,4/2,4/4,4/5,5/3,5/5}
    \DrawGivenMatrix{(0,0)}{5}{5}{\nnzcoef}
  }
}{Non-zero pattern of a matrix $A$}{matrix_pattern}{Example of non-zero pattern of a sparse matrix $A$}{matrix_and_pattern}
The sparse matrix $A\in\C^{n \times n}$ is stored in three one-dimensional arrays.
Two arrays of size $nnz$, $AA$ and $JA$, store respectively the non-zero values, through major row storage (row by row) and the column indices, thus $JA(k)_{1 \le k \le nnz}$ is the column index in $A$ matrix of $AA(k)_{1 \le k \le nnz}$.
Finally, $IA$, an array of size $n+1$ that stores the list of indices at which each row starts.
$IA(i)_{1 \le i \le n}$ and $IA(i+1)-1$ correspond respectively to the beginning and the end of the $i^{th}$ row in arrays $AA$ and $JA$, \ie, $IA(n+1) = nnz + 1$.
\includefig{
  \begin{eqnarray}
  AA & = & -5{\;}/{\;}14{\;}/{\;}8{\;}/{\;}1{\;}/{\;}2{\;}/{\;}10{\;}/{\;}4{\;}/{\;}2{\;}/{\;}9{\;}/{\;}15{\;}/{\;}7, \nonumber\\
  JA & = & 1{\;}/{\;}2{\;}/{\;}2{\;}/{\;}3{\;}/{\;}1{\;}/{\;}3{\;}/{\;}2{\;}/{\;}4{\;}/{\;}5{\;}/{\;}3{\;}/{\;}5, \nonumber\\
  IA & = & 1{\;}/{\;}3{\;}/{\;}5{\;}/{\;}7{\;}/{\;}10{\;}/{\;}12\nonumber
  \end{eqnarray}
}{Compressed-Sparse Row (\emph{CSR}) storage format}{}{csr_format}

  \subsection{Band-row splitting}

The partition of the equation set leads to allocating each processor a band of rows corresponding to the block of the processed vectors.
In \figreft{band-row_splitting} where an example of band-row splitting is given, these terms are located in a colorful area.
The band-row splitting approach consists in partitioning the matrix $A$ of size $n \times n$ into horizontal band matrices.
Each processor is in charge of the management of a band-row matrix of size $N_p \times n$ and the associated unknown vector $x$ of size $N_p \times 1$, as drawn in \figreft{band-row_splitting}.
This method of partitioning by band-row allows to exhibit a sufficient degree of properly balanced parallelism.
This implies assigning all processors, a block of equally sized rows, containing approximately the same number of non-zero coefficients.
Unfortunately, it suffers from a major lack of granularity for implementation on a distributed memory system.

\includefig{
  \SetDisplayIndex{1}
  \def\procSy{3}
  \def\procEy{6}
  \def\procSa{3}
  \def\procEa{6}
  \def\procSx{3}
  \def\procEx{6}
  \centering
  \setlength{\tabcolsep}{0.0pt}
  \begin{tabular}{ccc}
  $\qquad\qquad =$ & &\\
  \scalebox{0.55}{
  \SetSplitType{0}
  \def\nnzcoef{1/1,2/1,3/1,4/1,5/1,6/1,7/1,8/1,9/1,10/1}
  \SetDisplayInfoSize{1}
  \SetDisplayInfoSizeProc{1}
  \DrawGivenMatrixProcessor{(0,0)}{10}{1}{\nnzcoef}{\procSy}{\procEy}{\Large{$n$}}{\Large{$N_p$}}
  }
  &
  \scalebox{0.55}{
  \SetSplitType{0}
  \def\nnzcoef{
  1/1,1/2,1/4,1/6,
  2/1,2/2,2/3,2/5,
  3/2,3/3,3/4,3/5,
  4/1,4/3,4/4,4/7,4/10,
  5/3,5/5,5/7,
  6/1,6/6,6/7,6/9,
  7/3,7/5,7/6,7/7,7/10,
  8/2,8/4,8/8,8/10,
  9/6,9/9,
  10/4,10/7,10/8,10/10}
  \SetDisplayInfoSize{0}
  \SetDisplayInfoSizeProc{0}
  \DrawGivenMatrixProcessor{(0,0)}{10}{10}{\nnzcoef}{\procSa}{\procEa}{\Large{$n$}}{\Large{$N_p$}}
  }
  &
  \scalebox{0.55}{
  \SetSplitType{0}
  \def\nnzcoef{1/1,2/1,3/1,4/1,5/1,6/1,7/1,8/1,9/1,10/1}
  \SetDisplayInfoSize{1}
  \SetDisplayInfoSizeProc{1}
  \DrawGivenMatrixProcessor{(0,0)}{10}{1}{\nnzcoef}{\procSx}{\procEx}{\Large{$n$}}{\Large{$N_p$}}
  }\\
  $\;y$ & $\;\;\;\;\;\;\;A$ & $\;\;x$
  \end{tabular}
}{Example of band-row splitting of a matrix}{}{band-row_splitting}

  \subsubsection{Naive splitting}

The first splitting is the naive one, which consists in storing both matrices and vectors without any other informations.
The principle of the local matrix-vector product with row-band splitting consists in multiplying the local band row matrix by the global (\emph{gathered}) vector.
This implies data exchange between all processors (\emph{MPI\_ALLGATHER}).
To perform the matrix-vector product, each processor needs to receive and/or send missing information of the vector $x$ from/to cooperating processors.
According to the sparsity of the matrix, \ie, the distribution of non-zero values, the processor may receive and/or send unnecessary information which overload the communication.
This operation is very expensive for large-size matrices.
The communications dominate the computations for large-size matrices.
One solution to overcome this problem is to take into account the sparsity pattern of the local matrix, in order to send/receive only the necessary data from/to cooperating processors.\\

  \subsubsection{Sparsity Pattern splitting}

In this technique, in addition to the \emph{naive splitting} data, we store the information of dependencies from local to cooperating processors and from cooperating to local processors, according to the sparsity pattern of the local matrix.
Each processor has both a list of receive and a list of send dependencies, which keeps data exchange to a minimum.
\figreft{three_bandrow_splitting} gives an example of splitting a sparse matrix into three band rows.
\includefig{
  \SetGrid{0}
  \SetDisplayRowIndexD{0}{1}
  \SetDisplayColumnIndexD{1}{0}
  \def\NumbSubdom{3}
  \centering
  \setlength{\tabcolsep}{0.0pt}
  \begin{tabular}{cc}
    \scalebox{0.55}{
    \SetSplitType{0}
    \def\nnzcoef{
    1/1,1/2,1/4,1/6,
    2/1,2/2,2/3,2/5,
    3/2,3/3,3/4,3/5,
    4/1,4/3,4/4,4/7,4/10,
    5/3,5/5,5/7,
    6/1,6/6,6/7,6/9,
    7/3,7/5,7/6,7/7,7/10,
    8/2,8/4,8/8,8/10,
    9/6,9/9,
    10/4,10/7,10/8,10/10}
    \SetDisplayInfoSize{0}
    \SetDisplayInfoSizeProc{0}
    \DrawGivenMatrixProcessorSubdom{(0,0)}{10}{10}{\nnzcoef}{\NumbSubdom}
    }
    &
    \scalebox{0.55}{
    \SetSplitType{0}
    \def\nnzcoef{1/1,3/1,4/1,6/1,8/1,9/1,10/1}
    \SetDisplayInfoSize{1}
    \SetDisplayInfoSizeProc{1}
    \DrawGivenMatrixProcessorSubdom{(0,0)}{10}{1}{\nnzcoef}{\NumbSubdom}
    }\\
    $\;\;\;\;\;\;\;\;\;\;\;A$ & $\;\;\;\;\;\;x$
  \end{tabular}
}{Example of splitting of a matrix into three band rows}{}{three_bandrow_splitting}
\tablesubreft{list_dependency_sending} gives the corresponding list of dependencies for \textbf{sending} to cooperating processors of the splitting described in \figreft{three_bandrow_splitting}.
The corresponding list of dependencies for \textbf{receiving} from cooperating processors of the splitting described in \figreft{three_bandrow_splitting} is reported in \tablesubreft{list_dependency_receiving}.
\includesubtab{
{
  \renewcommand{\tabcolsep}{0.02cm}
  \begin{tabular}{|ccc|l|}
  \hline
  \textbf{Local proc.} &  & \textbf{Recv. proc} & \textbf{List Dependency nodes} \\
  \hline
   {\bfseries 1}       & $\leftarrow$ & 2          & 4 - \colorbox{gray}{5} - 6 \\
   \hline
   {\bfseries 2}       & $\leftarrow$ & 1          & 1 - 3 \\
   {\bfseries 2}       & $\leftarrow$ & 3          & \colorbox{gray}{7} - 9 - 10 \\
   \hline
   {\bfseries 3}       & $\leftarrow$ & 1          & 2 - 3 \\
   {\bfseries 3}       & $\leftarrow$ & 2          & 4 - \colorbox{gray}{5} - 6\\
   \hline
  \end{tabular}
}
}{\textbf{Receiving} dependencies}{list_dependency_receiving}{
{
  \renewcommand{\tabcolsep}{0.02cm}
  \begin{tabular}{|ccc|l|}
  \hline
  \textbf{Local proc.} & & \textbf{Send. proc} & \textbf{List Dependency nodes} \\
  \hline
  {\bfseries 1}        & $\rightarrow$ & 2          & 1 - 3 \\
  {\bfseries 1}        & $\rightarrow$ & 3          & 2 - 3 \\
  \hline
  {\bfseries 2}        & $\rightarrow$ & 1          & 4 - \colorbox{gray}{5} - 6\\
  {\bfseries 2}        & $\rightarrow$ & 3          & 4 - \colorbox{gray}{5} - 6 \\
  \hline
  {\bfseries 3}        & $\rightarrow$ & 2          & \colorbox{gray}{7} - 9 - 10\\
  \hline
  \end{tabular}
}
}{\textbf{Sending} dependencies}{list_dependency_sending}{List of \textbf{receiving} and \textbf{sending} dependencies of the splitting described in \figreft{three_bandrow_splitting}.}{list_dependency_receiving_sending}
The dependency nodes drawn in gray color in \tablereft{list_dependency_receiving_sending} correspond to zero values in vector $x$.
In fact, we can remove these nodes from the list of dependencies.
In this study, the band-row splitting with the ``\emph{Sparsity Pattern}'' technique takes into account both the sparsity of the matrix and the vector.
However, this technique does not guarantee a perfect load balance.
One solution for a perfect load-balance consists in partitioning the weighted graph (graph where the weight of the vertex $v_i$ associated with the row $i$ is its number of non-zero values) associated with the matrix in $k$ parts with the minimal amount of edges cut.
This solution leads to finding a good distribution of the sparse matrix on the parallel processors.
The concept consists in storing the $i^{th}$ row of the matrix on processor $j$ if the vertex $v_i$ is in the $i^{th}$ sub-part.
Then, we have equal weight in each band of the splitting.
The bands may be composed with non-contiguous rows.

  \paragraph{Matrix-vector product}

The processor that will perform the matrix-vector product for a band-row has only the corresponding terms of the vector $x$, the colored area in \figreft{band-row_splitting}.
In order to carry out the sparse matrix-vector, this process needs all the terms of the vector $x$. The first step, therefore, consists in collecting the terms that lacking from the colored area in \figreft{band-row_splitting}.
As it is the same for all processors, it will therefore be necessary to reconstruct the full vector $x$ on each processor.
This operation corresponds to a classic collective exchange, where each is both a transmitter and a receiver.
In this work, instead of using the collective operation, \emph{MPI\_Allgather}, including the message passing library (MPI)~\cite{gropp_usingmpiportable_2000}, we use the equivalent \emph{Send/Recv}, with a \emph{left-right} ordering of sending and receiving. For the processor $p$, the \emph{left-right} ordering consists in respectively sending and receiving to and from $k=p-1$, $k=p+1$, $k=p-2$, $k=p+2$, $k=p-3$, $k=p+3$, ..., if $k > 0$.  This process is described in \figimgreft{img/png/send_recv_ordering}.
\includefigimg[width=0.75\textwidth]{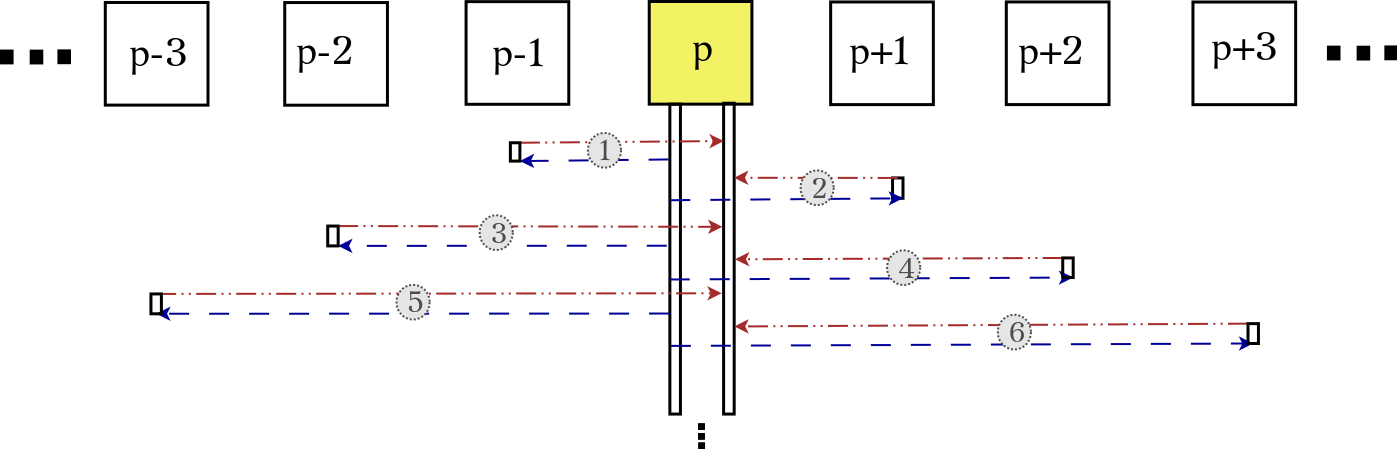}{Send/Recv ordering of the processor $p$}{}
The number of arithmetical operations requires to perform the local sparse matrix-vector multiplication, which is approximately $\frac{K \times n}{s}$, where $s$ is the number of processors, $n$ the dimension of the matrix, and $K$ the average number of non-zero coefficients per row.
On the other hand, the total number of terms of the vector $x$ to recover before performing the product is approximately $\frac{(s-1).n}{s}$, if the local matrix has non-zero values in almost all columns.
The amount of data is not small compared with the number of arithmetic operations.
Optimizing communications consists in finding a way to drastically limit the number of external values of vector $x$, located on the others processors, and is necessary to compute the product by the matrix.

  \paragraph{Basic linear algebra operations}

The computation of the dot product is a relatively simple operation. Each processor performs a local dot product, \ie, multiplies its elements and sums them, from their two local vectors.
Finally, the local sums are added using \emph{MPI\_ALLREDUCE} with the \emph{MPI\_SUM} operation.
Then each processor has the global dot product.
Operations such as an addition of vectors, the element wise product, etc. do not change compared to the sequential code.

  \subsection{Substructuring methods}
  \label{sec:Substructuring methods}

The substructuring method is the precursor of non-overlapping domain decomposition methods~\cite{letallec_domaindecompositionmethods_1994}~\cite{smith_domaindecompositionparallel_1996}~\cite{quarteroni_domaindecompositionmethods_1999}~\cite{toselli_domaindecompositionmethods_2005}.
The substructuring method is based on decomposition of the original structure into several sub-structures.
The term \textit{substructuring} is a way to describe the general method allowing to decompose splitting among subdomains sharing a common interface.
This method is most often used as a way to reduce the number of unknowns in the linear system by eliminating the interior unknowns.
.\\
Consider again the system of linear equations,
\begin{equation}
\label{eq:substructuring_linear_equation}
A u = f,
\end{equation}
where $A$ is a $n \times n$ square non-singular matrix, $f$ and $u$ represent respectively the right-hand side and the solution vector we are looking for.
Let us consider the field $\K$, usually $\R$ or $\C$.

  \subsubsection{Principle of substructuring methods}

In order to illustrate the substructuring method, we consider a problem steming from the finite element discretization of an elliptic partial differential problem.
To simplify the analysis, we consider the Laplace equation.
However, the analysis can be carried out for any coercive elliptic problem.
The model problem for the unknown $u$, in a bounded domain $\Omega$ with homogeneous Dirichlet boundary conditions on the boundary $\partial \Omega = \Gamma$ can be expressed as: for $f\in L^2(\Omega)$, find $u\in H^1(\Omega)$ such that $-\nabla^2 u = f\;in\; \Omega$ and $u = 0\;on\; \Gamma$.
An equivalent variational formulation of this problem can be formulated as: for $f\in L^2(\Omega)$, find $u\in H_0^1(\Omega)$ such that $\forall v \in H_0^1(\Omega), \int_\Omega \nabla u. \nabla v = \int_\Omega b.v$.
This problem is well posed, \ie it has one and only one solution.
After a Galerkin discretization with finite elements and a choice of nodal basis, the linear system~\eqref{eq:substructuring_linear_equation}, $Au = f$, is obtained, where $f$ denotes the right-hand side, $u$ the unknown and $A$ the stiffness matrix which is sparse.
The following shows the main steps for performing a matrix-vector product using a splitting into sub-structures.\\

\subsubsection{Matrix splitting}

In practice, mesh partitioning is a crucial step of the finite element method.
A finite element matrix is associated with a finite element mesh and the elements of the matrix are correlated with the interaction of the basis functions defined in the elements of the mesh. The total matrix is calculated as an assembly of elementary matrices.
Let us consider a global domain $\Omega$ partitioned into two subdomains without overlap $\Omega_1$ and $\Omega_2$, with a shared interface $\Gamma$ as drawn in \figimgreft{img/png/substructuring_domain}.
\includefigimg[width=0.55\textwidth]{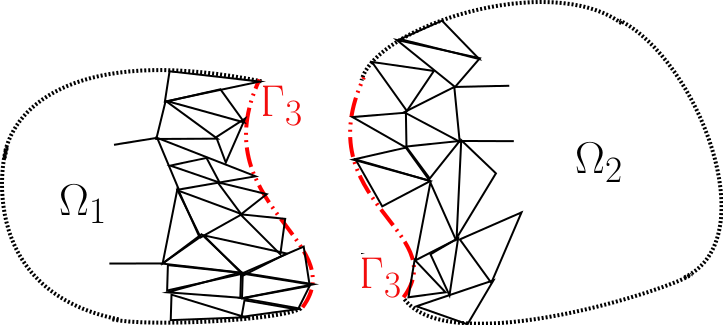}{Splitting into two subdomains}{}
When a suitable numbering of the degrees of freedom is harnessed, the stiffness matrix of the initially considered model problem can be written as in the following matrix:
\begin{eqnarray}
\label{eqn:2.4:eqn0}
  A =
  \left( \begin{array}{ccc}
  \Aij{A}{1}{1} & 0      & \Aij{A}{1}{3} \\
  0      & \Aij{A}{2}{2} & \Aij{A}{2}{3} \\
  \Aij{A}{3}{1} & \Aij{A}{3}{2} & \Aij{A}{3}{3}
  \end{array} \right)
\end{eqnarray}
It is formulated considering the case where the set of nodes numbered $1$ and $2$ are respectively associated to the subdomains $\Omega_1$ and $\Omega_2$.
The last set of nodes numbered $3$ corresponds to the interface nodes of both subdomains.
It is possible to split the original system~\eqref{eq:substructuring_linear_equation} into blocks
\begin{eqnarray}
\label{eqn:2.4:eqn1}
  \left( \begin{array}{ccc}
  \Aij{A}{1}{1} & 0      & \Aij{A}{1}{3} \\
  0      & \Aij{A}{2}{2} & \Aij{A}{2}{3} \\
  \Aij{A}{3}{1} & \Aij{A}{3}{2} & \Aij{A}{3}{3}
  \end{array} \right)
  \left( \begin{array}{ccc}
  \ri{u}{1} \\
  \ri{u}{2} \\
  \ri{u}{3}
  \end{array} \right)
  =
  \left( \begin{array}{ccc}
  \ri{f}{1} \\
  \ri{f}{2} \\
  \ri{f}{3}
  \end{array} \right)
\end{eqnarray}
where $u=(\ri{u}{1},\ri{u}{2},\ri{u}{3})^t$ is the unknown vector and $f=(\ri{f}{1},\ri{f}{2},\ri{f}{3})^t$ is the right-hand side.
The blocks $\Aij{A}{1}{3}$ and $\Aij{A}{2}{3}$ are respectively the transpose matrix of $\Aij{A}{3}{1}$ and $\Aij{A}{3}{2}$, and the blocks $\Aij{A}{1}{1}$ and $\Aij{A}{2}{2}$ are symmetric positive-definite if $A$ was symmetric positive-definite.
By assigning the different subdomains at distinct processors, the local matrices can be formulated in parallel as follows:
\begin{eqnarray}
\label{eqn:2.4:mat_subdomains}
  A_{1} =
  \left( \begin{array}{cc}
  \Aij{A}{1}{1} & \Aij{A}{1}{3} \\
  \Aij{A}{3}{1} & \Aijs{A}{3}{3}{1}
  \end{array} \right)
  \quad\text{ and }\quad
  A_{2} =
  \left( \begin{array}{cc}
  \Aij{A}{2}{2} & \Aij{A}{2}{3} \\
  \Aij{A}{3}{2} & \Aijs{A}{3}{3}{2}
  \end{array} \right).
\end{eqnarray}
The blocks $\Aijs{A}{3}{3}{1}$ and $\Aijs{A}{3}{3}{2}$ denote the interaction between the nodes on the interface $\Gamma$, respectively integrated in subdomains $\Omega_1$ and on $\Omega_2$, \ie
\begin{equation}
\label{eq:interface_summation}
\Aij{A}{3}{3} = \Aijs{A}{3}{3}{1} + \Aijs{A}{3}{3}{2}
\end{equation}
In practice, the subdomains $\Omega_1$ and $\Omega_2$ respectively know the set of nodes ($1$, $3$) and ($2$,$3$).

  \subsubsection{Matrix-vector product}

As described in~\cite{cheikahamed_fastsparsematrix_2012}~\cite{cheikahamed_iterativemethodssparse_2012}~\cite{magoules_fastiterativesolvers_2015}, iterative Krylov algorithms require performing one or more multiplication(s) of the matrix $A$ by a descent direction vector $u=(\ri{u}{1},\ri{u}{2},\ri{u}{3})^t$ at each iteration.
With the splitting into two subdomains, the global matrix-vector multiplication can be written as follows:
\begin{eqnarray}
\label{eqn:2.4:spmv}
  \left( \begin{array}{ccc}
  \ri{y}{1} \\
  \ri{y}{2} \\
  \ri{y}{3}
  \end{array} \right)
  & = &
  \left( \begin{array}{ccc}
  \Aij{A}{1}{1} & 0      & \Aij{A}{1}{3} \\
  0      & \Aij{A}{2}{2} & \Aij{A}{2}{3} \\
  \Aij{A}{3}{1} & \Aij{A}{3}{2} & \Aij{A}{3}{3}
  \end{array} \right)
  \left( \begin{array}{ccc}
  \ri{u}{1} \\
  \ri{u}{2} \\
  \ri{u}{3}
  \end{array} \right) \nonumber\\
  & = &
  \left( \begin{array}{ccc}
  \Aij{A}{1}{1}\ri{u}{1} + \Aij{A}{1}{3}\ri{u}{3} \\
  \Aij{A}{2}{2}\ri{u}{2} + \Aij{A}{2}{3}\ri{u}{3} \\
  \Aij{A}{3}{1}\ri{u}{1} + \Aij{A}{3}{2}\ri{u}{2} + \Aij{A}{3}{3}\ri{u}{3}
  \end{array} \right)\nonumber
\end{eqnarray}
Considering the local matrices' described in equation~(\ref{eqn:2.4:mat_subdomains}), we can independently compute both local matrix-vector products as follows
\begin{eqnarray}
  \left( \begin{array}{c}
  \ri{y}{1}\\
  \ris{y}{3}{1}
  \end{array} \right)
  =
  \left( \begin{array}{cc}
  \Aij{A}{1}{1}\ri{u}{1} + \Aij{A}{1}{3}\ri{u}{3} \\
  \Aij{A}{3}{1}\ri{u}{1} + \Aijs{A}{3}{3}{1}\ri{u}{3}
  \end{array} \right)\nonumber
\end{eqnarray}
\begin{eqnarray}
  \left( \begin{array}{c}
  \ri{y}{2}\\
  \ris{y}{3}{2}
  \end{array} \right)
  =
  \left( \begin{array}{cc}
  \Aij{A}{2}{2}\ri{u}{2} + \Aij{A}{2}{3}\ri{u}{3} \\
  \Aij{A}{3}{2}\ri{u}{2} + \Aijs{A}{3}{3}{2}\ri{u}{3}
  \end{array} \right)\nonumber
\end{eqnarray}
Since $\Aij{A}{3}{3} = \Aijs{A}{3}{3}{1} + \Aijs{A}{3}{3}{2}$, and $\ri{y}{3} = \ris{y}{3}{1} + \ris{y}{3}{2}$. According to this last remark, SpMV can be calculated in two steps:
\begin{itemize}
  \item calculate the local matrix-vector multiplication in each subdomain
  \item assemble, the local contributions on the interface
\end{itemize}
The first step involves only local data.
The second requires the exchange of data between processes dealing with subdomains with a common interface.
In order to assemble interface values of neighboring subdomains, each processor responsible for a subdomain must know the description of its interfaces.

\subsubsection{Exchange at the interfaces}

When a subdomain $\Omega_i$ has several neighboring subdomains, we denote $\Gamma_{ij}$ the interface between $\Omega_i$ and $\Omega_j$ as described in \figimgref{img/png/substructuring_interface}.
\includefigimg[width=0.35\textwidth]{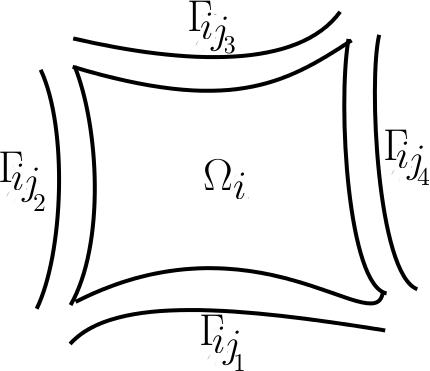}{Sub-structure interface description}{}
An interface is identified by its neighboring subdomains and the equations associated with its nodes.
The interface is evaluated from its subdomains using the sparse matrix-vector product.
This computation is in two steps for each neighboring subdomain (\algorithmreft{collect_send_neighbor}): \underline{\emph{collect}} the values of the local vector $y = Au$ for all interfaces nodes, and then \underline{\emph{send}} this list to the $y$ vector of the interface equation.
\includealgo[
  \Input{$\;n_s$: number of interface nodes\\
         $\;number\_of\_neighboring$: number of neighboring subdomains\\
         $\;y$: values on the whole subdomain\\
         $\;list_s$: list of the interface nodes}
  \Output{$\;buffer_s$: sending buffer}
  \algovariable $s$, $i$\\
]{
  \For{$s \gets 1:number\_of\_neighboring$}{
    \For{$i = 1:n_s$}{
      $buffer_s(i) = y(list_s(i))$
    }
    Send $buffer_s$ to neighbour(s)
  }
}{Construct inner buffer and send to neighboring}{}{collect_send_neighbor}
The next step consists in updating these changes to all neighboring subdomains at interfaces equations.
First, the contributions of the array containing the result of the matrix-vector product at the interface are received, and then values on the corresponding interface nodes are updated.
This processus is described in \algorithmreft{recv_update}.
\includealgo[
  \Input{$\;n_s$: number of interface nodes\\
         $\;number\_of\_neighboring$: number of neighboring subdomains\\
         $\;y$: values on the whole subdomain\\
         $\;list_s$: list of the interface nodes}
  \Output{$\;buffer_s$: receiving buffer}
  \algovariable $s$, $i$\\
]{
  \For{$s \gets 1:number\_of\_neighboring$}{
    Receive $buffer_s$ from neighbour(s)
    \For{$i = 1:n_s$}{
      $y(list_s(i)) = y(list_s(i)) + buffer_s(i)$
    }
  }
}{Receiving interface results and updating interface equations}{}{recv_update}
When an equation is shared by several interfaces, the node value of the local vector $y = Au$  in question is sent to all interfaces to which it belongs.
For any number of subdomains, the mechanism of interface exchange and update is similar to those previously presented.

The use of the substructuring approach in iterative algorithms is inherently parallel and makes it an excellent candidate for implementation on parallel computers.
Indeed, we can distribute the subdomains over all available processors and thus compute the matrix-vector products locally, independently and in parallel, and use the distributed memory in order to limit the memory usage.
As explained previously, after the computing of the local matrix-vector multiplications, they required to be assembled along the interface.
The key ingredient of the data is the local matrix $C$ that arises from the finite element discretization. With this approach, each node $i$ only needs to store $C_i$, the corresponding local matrix to the subdomain $\Omega_i$, which is only a fraction of the original matrix.
The dot product requires that each processor compute a weighted combination of the interface contributions in order to update its own data.
After that, an \emph{MPI\_ALLREDUCE} is required to compute the global inner product.

The iterative substructuring method introduces only two new steps, which reside in data exchange. They consist firstly in sharing the contributions of the local computed SpMV at the interface.
Each machine requires to know the list of nodes along the interface and the number of neighboring subdomains.
Secondly, in assembling results over the cluster in order to piece together the local scalar product.
This action, realized with MPI, is independent from the splitting.
Finally, another advantage of this algorithm is that it can easily be generalized for $n$ subdomains.
This approach however, presents two disadvantages.
The first drawback arises from a computational point of view.
The granularity, \ie, the number of operations to be performed by the processors compared to the amount of data received or sent by the processors may be weak.
Indeed, here the granularity is proportional to the number of nodes in the subdomains compared to the number of nodes on the interface.
The number of operations depend on the first parameter and the data transfer depends on the second parameter.
If a lot of subdomains are used, the interface size will not be small compared to the local sub-problem size.
This means that the processors realize few computing operations (a local matrix-vector product) and a lot of communications.
The second and more important drawback is an algorithmic one.
The classical parallel preconditioners per subdomain are based on an incomplete factorization of the local matrices.
Such preconditioners are less and less efficient when the number of subdomain increases.

\section{Benchmark environments}
\label{sec:Benchmark environments}

  \subsection{Alinea: an hybrid CPU/GPU library}

In this section, we present the main features of the library, Alinea, we have implemented and used in this paper to perform effectivly advanced linear algebra and solving linear systems on both CPU and GPU.
Alinea stands for \emph{{\bf A}dvanced {\bf LINE}ar {\bf A}lgebra}.
Alinea is targeted as a scalable software for proposing effective linear algebra operations on both CPU and GPU platforms.
It includes numerous algorithms for solving linear systems, with different matrix storage formats, with real and complex arithmetics in single and double precision, on both CPU and GPU devices.
Alinea is devoted to simplifying the development of engineering and science problems on CPU and GPU by discharging most of the difficulties encountered when using these architectures, particularly with GPU.
Alinea investigates and seeks the best way to effectively implement linear algebra operations and solver algorithms on both CPU and GPU. It also allows to write the same code for the CPU and GPU versions of an algorithm easily.
The library is implemented in C++ language and proposes a simple C API interface.
The main features of Alinea described in \tabletblreft{Alinea_main_level} which are applied to real and complex arithmetic numbers, simple and double precision.
\includetabtbl{
  {
    \begin{tabular}{cccc}
      \hline
      type          & function                   & example        & reference\\
      \hline
      vector-vector & addition \& multiplication & $w=u+\alpha v$ & BLAS 1\\
      vector-vector & scalar product             & $p=<u,v>$      & BLAS 1\\
      vector        & norm                       &  $q=||u||$     & BLAS 1\\
      matrix-vector & sparse product             & $y=A*x$        & BLAS 2\\
      matrix-matrix & matrix product             & $C=A*B$        & BLAS 3\\
      solver        & direct                     & $x=LU(A,b)$    & \\
      solver        & iterative                  & $x=CG(A,b)$    & \\
      \hline
    \end{tabular}
  }
}{Alinea (\emph{{\bf A}dvanced {\bf LINE}ar {\bf A}lgebra}) main levels}{}{Alinea_main_level}
We are dealing with huge matrix systems and therefore parallel computing becomes an important issue for optimization.
The GPU version with CUDA or OpenCL, written in C++, has the task to load the methods into the GPU, to control the data transfer and to manage the memory.
This concept allows to develop a much more modular code, easier to use for users or developers.
Indeed, each block is independent and the use of templates allows a very intuitive use of the library, regardless of the architecture or the data type selected (real, complex, simple precision or double-precision).
Therefore, the available C++ methods are identical.
The template is designed to \lstinline+<T,U>+, where \lstinline+T+ is the type of value, \eg, \lstinline+double+ or \lstinline+std::complex<double>+ and \lstinline+U+ the type of index, \eg, \lstinline+int+ or \lstinline+unsigned int+.

  \subsection{Hardware configuration}

The performance of our algorithms depends on the characteristics of the test machines~\cite{evans_parallellinearsystem_1979}~\cite{cheikahamed_fastsparsematrix_2012}~\cite{margaris_parallelimplementationsjacobi_2014}~\cite{magoules_alineaadvancedlinear_2015}. The calculations depend on the precision (single or double) and the hardware support.
In this paper, the experiments have been performed on \emph{MRG/LISA} cluster, which is a in-house cluster of our team.
\emph{MRG/LISA} is an hybrid CPU/GPU cluster, whose characteristics are given in \tabletblreft{mrg_lisa_feature}.
The cluster consists of 6 CPU nodes $(C_1,C_2,\ldots,C_6)$ and 8 CPU/GPU nodes $(G_1,G_2,\ldots,G_8)$.
The version of the MPI library used is OpenMPI \emph{(OpenRTE) 1.6.5}.
\includetabtbl{
  {
    \tiny
    \renewcommand{\tabcolsep}{0.05cm}
    \begin{tabular}{|clll|}
      \hline
      \textbf{node}  & \textbf{OS} & \textbf{CPU} & \textbf{GPU}\\
      \hline
      \hline
      \begin{tabular}{c}
      $\bs{C_1,C_2,C_3}$\\\\
      $\bs{C_4,C_5,C_6}$
      \end{tabular}
      & 
      \begin{tabular}{l}
      Linux $64$ bits\\
      Ubuntu 14.04 LTS\\
     \emph{system '/'}: $55GB$ \\
     \emph{swap}: $30GB$ \\
     \emph{/home}: $180GB$
      \end{tabular}
      & 
      \begin{tabular}{l}
      Intel(R) Xeon(R) E5410\\
      2,33GHz, $2 \times 4 = 8$ cores \\
      RAM: $8$ GB\\
      MPI (OpenRTE) 1.6.5
      \end{tabular}
      & 
      none\\
      \hline
      \begin{tabular}{c}
      $\bs{G_1,G_2,G_3,G_4}$
      \end{tabular}
      & 
      \begin{tabular}{l}
      Linux $64$ bits\\
      Ubuntu 14.04 LTS\\
     \emph{system '/'}: $55GB$ \\
     \emph{swap}: $30GB$ \\
     \emph{/home}: $180GB$
      \end{tabular}
      & % --
      \begin{tabular}{l}
      Intel(R) Core(TM) i7 \\
      2,80GHz, $2 \times 4 = 8$ cores \\
      RAM: $8$ GB\\
      MPI (OpenRTE) 1.6.5
      \end{tabular}
      & % --
     \begin{tabular}{l}
        Tesla K20c $4799MB$ \\
        GTX 570 $1279MB$ \\
        CUDA v6.5 \\
        Double precision\\
     \end{tabular}\\
      \hline
       % ------------------------------------------------------------------------
      \begin{tabular}{c}
      $\bs{G_5,G_6,G_7,G_8}$
      \end{tabular}
      & % --
      \begin{tabular}{l}
      Linux $64$ bits\\
      Ubuntu 14.04 LTS\\
     \emph{system '/'}: $55GB$ \\
     \emph{swap}: $30GB$ \\
     \emph{/home}: $180GB$
      \end{tabular}
      & % --
      \begin{tabular}{l}
      Intel(R) Xeon(R) E5-2609\\
      2,10GHz, $4 \times 6 = 24$ cores \\
      RAM: $16$ GB\\
      MPI (OpenRTE) 1.6.5
      \end{tabular}
      & % --
     \begin{tabular}{l}
        $\bs{G_5,G_6,G_7}$\\
        Quadro K4000 $3071 MB$ \\
        $\bs{G_8}$\\
        Quadro K600 $1023MB$\\
        $\bs{G_5,G_6,G_7,G_8}$\\
        CUDA v6.5 \\
        Double precision\\
     \end{tabular}\\
      \hline
      \hline
       \multicolumn{4}{|c|}{The interconnected network is a switched, star shaped 10Mb/s Ethernet network.}\\
      \hline
    \end{tabular}
  }
}{MRG/LISA hybrid computing clusters}{}{mrg_lisa_feature}

\section{Numerical results}
\label{sec:Numerical results}

In this section, we report the numerical experiments of parallel synchronous Jacobi method for different splitting models include naive band-row splitting (BANDROW (JB)), optmized band-row splitting (BANDROW-OP (JBO)) and substructuring methods (SUBSTRUCTURING (JSS)).
The execution times reported in this part correspond to the average time of 10 executions.

  \subsection{Matrices tested}

To evaluate the comparison of the three splitting strategies presented in this paper, we consider the 3D Laplace equation (Poisson equation):
\begin{equation}
  \label{eq:luf_laplace}
  - \Delta u = f,
\end{equation}
in the domain $\Omega = [0; 250,000] \times [0; 250,000] \times [0; 250,000]$, which is projected into $[0; 1] \times [0; 1] \times [0; 1]$.
The projection conserves the geometry, the discretization topology ($h$ remains uniform), the number of nodes and elements, the layers and the partitioning.

  \subsubsection{Set-1: 3D Laplace equation, academic test cases}

The first set of data consists of academic test cases, which consider that $u(x,y,z)$ satisfies the following Dirichlet boundary conditions,
\begin{equation}
\label{eq:luf_laplace_dirichlet}
\begin{cases}
u(x,y,z) = 0, \quad \forall\ (x,y,z) \in \{0;1\} \times [0; 1] \times [0; 1],\\
u(x,y,z) = 0, \quad \forall\ (x,y,z) \in [0; 1] \times \{0;1\}  \times [0; 1],\\
u(x,y,z) = 0, \quad \forall\ (x,y,z) \in [0; 1] \times [0; 1] \times \{0;1\}.
\end{cases}
\end{equation}
The right hand side is defined such that $f(x,y,z) = \cos(x + y)$.
\includetabtbl{
  {
  \begin{tabular}{l||lr}
    \hline
    \textbf{}  & \textsc{cube-35937} & \textsc{cube-274624} \\
    \hline
    instance & \multicolumn{1}{c}{\figsubimgreft{img_cube-35937--00004}} & \multicolumn{1}{c}{\figsubimgreft{img_cube-274624--00004}}\\
    $x-$range & $[0:0.03125:1]$ & $[0:0.015625:1]$\\
    $y-$range & $[0:0.03125:1]$ & $[0:0.015625:1]$\\
    $z-$range & $[0:0.03125:1]$ & $[0:0.015625:1]$\\
    Number of nodes & $35,937$ & $274,624$\\
    \hline
  \end{tabular}
  }
}{Statistics of the academic cube meshes: \emph{luf\_cube-35937} and \emph{luf\_cube-274624}}{}{statistics_cube_meshes}
Thist first set includes two matrices of size $35,937$ and $274,624$ obtained from the finite element discretization of the equation~\eqref{eq:luf_laplace}.
The features of their corresponding meshes are reported in \tabletblreft{statistics_cube_meshes}.
\includesubfigimg{scale=0.25}{img_cube-35937--00004}{luf\_cube-35937}{scale=0.25}{img_cube-274624--00004}{luf\_cube-274624}{Finite element mesh examples of the cube}{mesh_cube-35937-274624--00004}
\figimgreft{mesh_cube-35937-274624--00004} gives an example of the finite element mesh of the \textsc{cube-35937} \figsubimgref{img_cube-35937--00004} and \textsc{cube-274624} \figsubimgref{img_cube-274624--00004}.

  \subsubsection{Set-2: 3D gravitational potential equation}

The second set of data are obtained from the finite element discretization of the gravitational potential equation.
The gravitational potential of a density anomaly distribution is a particular case of the equation~\eqref{eq:luf_laplace} with a right hand side defined such that $f(x,y,z) =4\pi G \delta\rho(x,y,z)$, where $\delta\rho$ is the density anomaly and $G$ the gravitational constant.
The domain $\Omega = [0; 250,000] \times [0; 250,000] \times [-15,000; 0]$ is related to the region of Chicxulub impact crater, which is located underneath the town of Chicxulub in Yucat{\'a}n, southwest of Mexico on the Yucat{\'a}n Peninsula.
The domain $\Omega$ is projected into $[0; 1] \times [0; 1] \times [0; 0.06]$.

The features of this domain are presented by the authors in~\cite{magoules_optimizedschwarzmethod_2016}.
We also consider homogeneous Dirichlet boundary conditions, which are defined on the vertical faces ($xz-$plane and $yz-$ plane), \ie
\begin{equation}
\label{eq:luf_gravitational_dirichlet}
\begin{cases}
u(x,y,z) = 0, \quad \forall\ (x,y,z) \in \{0;1\} \times [0; 1] \times [0, 0.06],\\
u(x,y,z) = 0, \quad \forall\ (x,y,z) \in [0; 1] \times \{0;1\}  \times [0, 0.06]
\end{cases}
\end{equation}
The density anomaly $\delta\rho$ are coming from geological surveys of Chicxulub impact crater.
The mesh characteristics of the \emph{Set-2} used for the finite element discretization of the gravitational potential equation are reported in~\tablereft{statistics_original_projected_meshes}.
\includesubtab{
{
  \tiny
  \renewcommand{\tabcolsep}{0.12cm}
  \begin{tabular}{l||lr}
    \hline
    \textbf{}  & \multicolumn{2}{c}{\textbf{original mesh}}\\
    \textbf{}  & \textsc{level \#0} & \textsc{level \#1} \\
    \hline
    $x-$range & $[0:2500:250,000]$ & $[0:1250:250,000]$ \\
    $y-$range & $[0:2500:250,000]$ & $[0:1250:250,000]$ \\
    $z-$range & $[-15,000:2500:0]$ & $[-15,000:1250:0]$ \\
    Number of nodes & $71,407$ & $525,213$ \\
    Number of elements & $60,000$ & $480,000$ \\
    Size of element & $\approx 4330.12702$ & $\approx 2165.06351$ \\
    Total interface size & $96, 234$ & $391,510$ \\
    \hline
  \end{tabular}
}
}{\textbf{Original} meshes}{statistics_original_meshes}{
{
  \tiny
  \renewcommand{\tabcolsep}{0.12cm}
  \begin{tabular}{l||lr}
    \hline
    \textbf{}  & \multicolumn{2}{c}{\textbf{projected mesh}}\\
    \textbf{}  & \textsc{level \#0} & \textsc{level \#1} \\
    \hline
    $x-$range & $[0:0.01:1]$ & $[0:0.005:1]$ \\
    $y-$range & $[0:0.01:1]$ & $[0:0.005:1]$ \\
    $z-$range & $[0:0.01:0.06]$ & $[0:0.005:0.06]$ \\
    Number of nodes & $71,407$ & $525,213$\\
    Number of elements & $60,000$ & $480,000$\\
    Size of element & $\approx 0.0173$ & $\approx 0.0087$\\
    Total interface size & $96, 234$ & $391,510$\\
    \hline
  \end{tabular}
}
}{\textbf{Projected} meshes}{statistics_projected_meshes}{Statistics of the original and projected meshes}{statistics_original_projected_meshes}
The matrices are identified by \emph{gravi\_hexas100x100x6\_0} (level \#1) and \emph{gravi\_hexas100x100x6\_1} (level \#2).

  \subsubsection{Sketches of matrices of test cases}

\tabletblreft{sketches_matrices_fem_test_cases} sums the tested matrices of the numerical expermiments.
The table reports the main properties of each matrix: $h$ the size of the matrix, $nnz$ the number of non-zero values, $density$ the density that corresponds to the number of non-zero values divided by the total number of matrix coefficients, $nnz/h$ the mean row density, $max\_row$ the maximal row density,  $\sigma(nnz/n)$ or $nnz/h\;stddev$ the standard deviation of the mean row density ($nnz/h$) and $bandwidth$ the upper bandwidth, which is equal to the lower bandwidth in the case of a symmetric matrix.
The first and second pictures represent respectively for each matrix the pattern of non-zero values and an histogram of the distribution of non-zero values per row.
\newcommand{\infosparseT}[9]{
\includegraphics[scale=0.11]{sparse/#1_sparse}\ \includegraphics[scale=.10]{sparse/#1_sNnzRow}
& \begin{tabular}{m{2.0cm}m{2.5cm}}
\multicolumn{2}{c}{\textbf{#2}}\\
h = #3 & density = #5\\
nz = #4 & bandwidth = #6\\
max row = #7 & nz/h = #8\\
nz/h stddev = #9 & \\
\end{tabular}
}
\includetabtbl{
  {
  \renewcommand{\tabcolsep}{0.03cm}
\begin{tabular}{m{3.4cm}m{3.5cm}}
\hlinewd{1.0pt}
\infosparseT{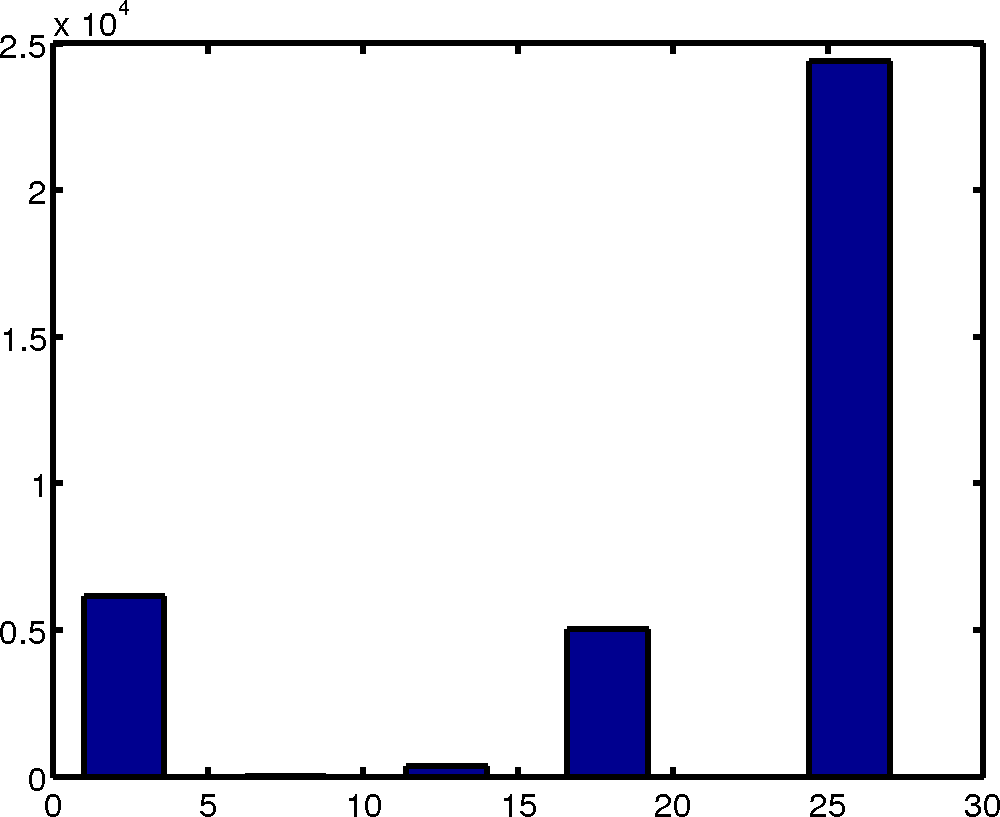}{luf\_cube-35937}{35,937}{759,667}{0.059}{1,123}{27}{21.139}{9.741}\\
\multicolumn{2}{m{8.3cm}}{\emph{Structured FEM problem, 3D Laplace equation, cube $[0:0.03125:1]\times[0:0.03125:1]\times[0:0.03125:1]$.}}\\
\hlinewd{1.0pt}
\infosparseT{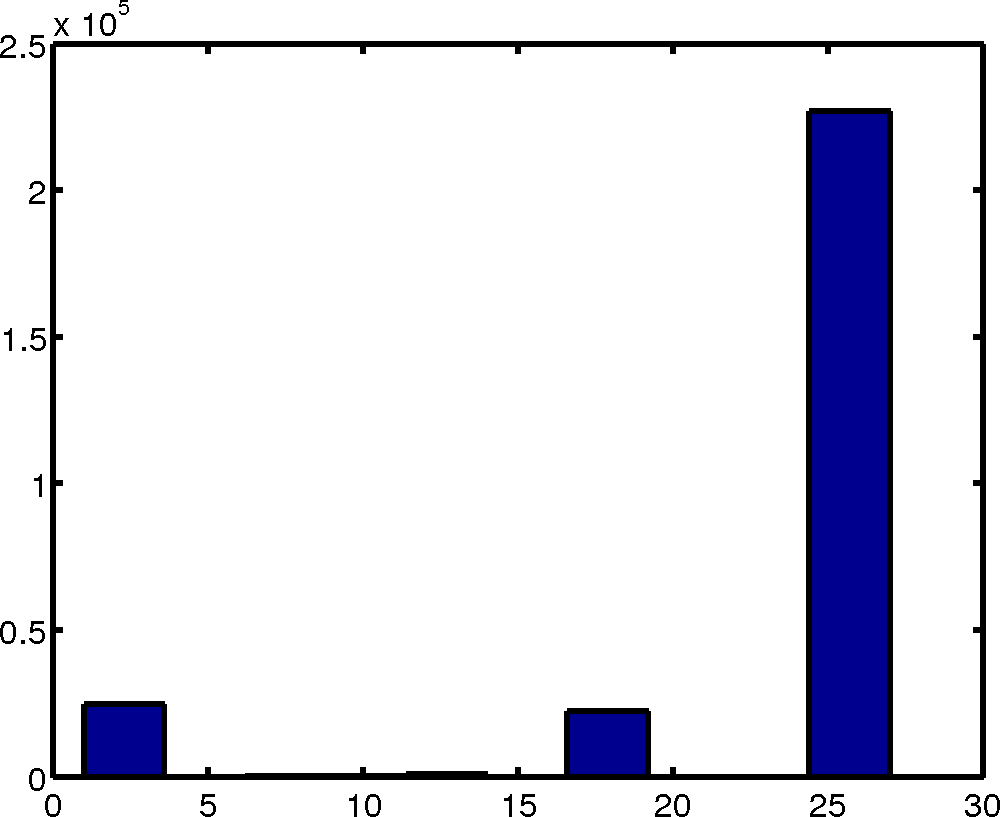}{luf\_cube-274625}{274,625}{6,563,167}{0.009}{4,291}{27}{23.899}{7.621}\\
\multicolumn{2}{m{8.3cm}}{\emph{Structured FEM problem, 3D Laplace equation, cube $[0:0.015625:1]\times[0:0.015625:1]\times[0:0.015625:1]$.}}\\
\hlinewd{1.0pt}
\infosparseT{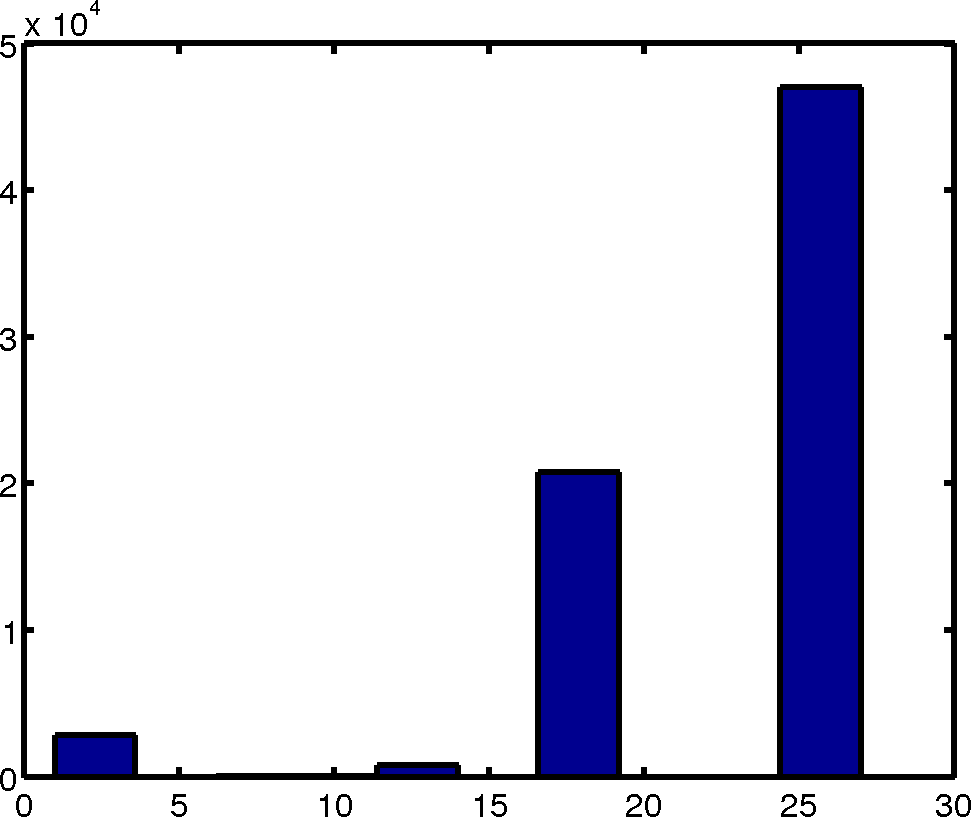}{gravi\_hexas100x100x6\_0}{71,407}{1,656,131}{0.032}{10,303}{27}{23.193}{6.174}\\
\multicolumn{2}{m{8.3cm}}{\emph{Structured FEM problem, 3D gravitational potential equation, parallelepiped $[0:0.01:1]\times[0:0.01:1]\times[0:0.01:0.06]$.}}\\
\hlinewd{1.0pt}
\infosparseT{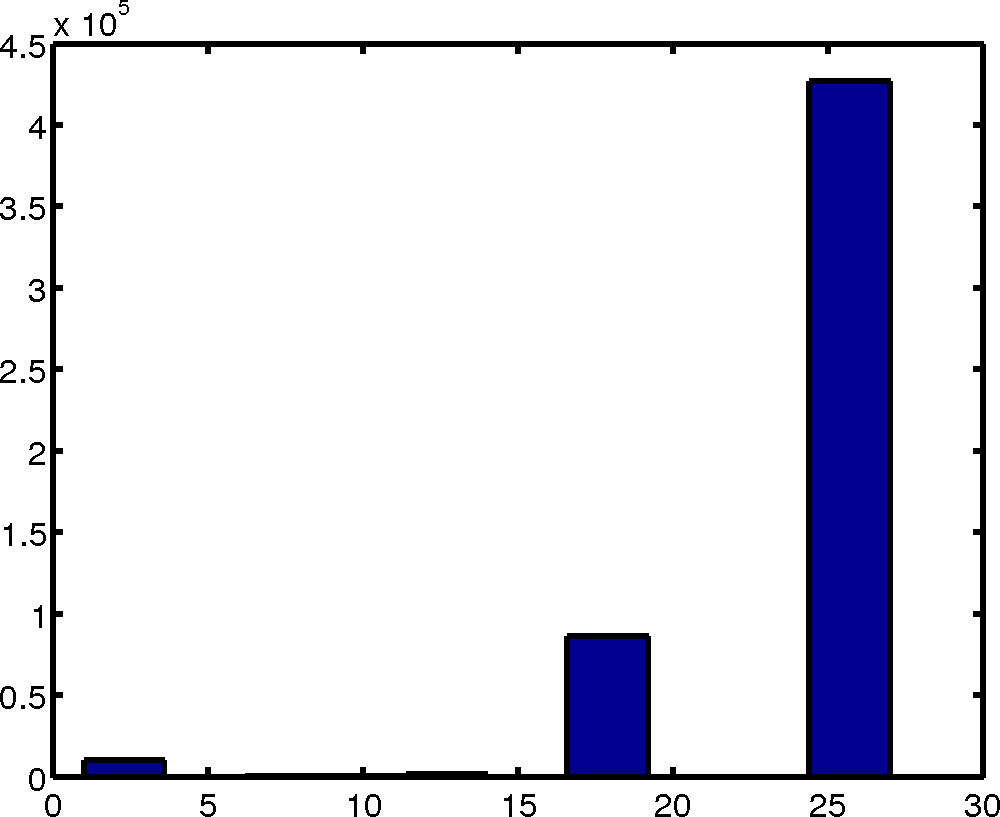}{gravi\_hexas100x100x6\_1}{525,213}{13,107,979}{0.005}{463,803}{27}{24.957}{4.820}\\
\multicolumn{2}{m{8.3cm}}{\emph{Structured FEM problem, 3D gravitational potential equation, parallelepiped $[0:0.005:1]\times[0:0.005:1]\times[0:0.005:0.06]$.}}\\
\hlinewd{1.0pt}
\end{tabular}
  }
}{Sketches of matrices obtained with the finite element discretization of the \emph{Set-1} and \emph{Set-2}}{}{sketches_matrices_fem_test_cases}

  \subsection{Experimental results}

The system is solved in parallel using a Jacobi splitting (Section~\ref{sec:Jacobi iterations}).
The benchmark consists in an analysis of the Jacobi method declined in 3 versions: \emph{BANDROW (JB)}, \emph{BANDROW-OP (JBO)} and \emph{SUBSTRUCTURING (JSS)}, which correspond respectively to the Jacobi method with naive band-row partitioning, optimized band-row partitioning and substructuring partitioning (Section~\ref{sec:Matrix partitioning}).
The experiments have been performed on \emph{LISA cluster}, which is described in~\tabletblreft{mrg_lisa_feature} with the set of matrices described in~\tabletblreft{sketches_matrices_fem_test_cases}. We have used the nodes $(C_1,C_2,C_3,C_4,C_5,C_6)$.

The convergence criterion is the weighted norm defined by $\norm[\infty]{D\left( \uk{k+1} - \uk{k}\right) } < \mathepsilon$, where the residual threshold $\mathepsilon = 10^{-8}$.
The norm is weighted to the diagonal $D$.

\tabletblreft{parallel_luf_cube-35937}, \tabletblreft{parallel_luf_cube-274625} and \tabletblreft{parallel_gravi_hexas100x100x6_0}, \tabletblreft{parallel_gravi_hexas100x100x6_1} show the numerical results of the experiments.
These tables report respectively the experimental results of the \emph{BANDROW (JB)}, \emph{BANDROW-OP (JBO)} and \emph{SUBSTRUCTURING (JSS)} synchronous Jacobi method.

These tables are organized as follows.
The number of processors is given in the first column, and the results of the \emph{BANDROW (JB)}, \emph{BANDROW-OP (JBO)} and \emph{SUBSTRUCTURING (JSS)} synchronous Jacobi method are respectively reported from column 2 to 4, from column 5 to 7, and from column 8 to 10.
For each method, the first column gives the number of iterations (\# iter), the second column collects respectively the communication (exchange) times and the solver execution times in seconds (s). The third column gives the efficiency (eff) of the synchronous algorithm upon the sequential code (one processor).
\includetabtbl{
  {
  \renewcommand{\arraystretch}{1.0}
  \renewcommand{\tabcolsep}{0.05cm}
  \begin{tabular}{|>\bfseries c|ccc|ccc|ccc|}
    \hlinewd{1.0pt}
    & \multicolumn{9}{c|}{\textbf{\textsc{Synchronous}}} \\
    & \multicolumn{3}{c|}{\emph{J-BANDROW}} & \multicolumn{3}{c|}{\emph{J-BANDROW-OP}} & \multicolumn{3}{c|}{\emph{J-SUBSTRUCTURING}}\\
    \hline
    \textbf{\#p (\#n)} & {\bfseries \# iter.} & {\bfseries time (s)} & {\bfseries eff.} & {\bfseries \# iter.} & {\bfseries time (s)} & {\bfseries eff.} & {\bfseries \# iter.} & {\bfseries time (s)} & {\bfseries eff.} \\
                       & \psync{}{comm}{total}{}& \psync{}{comm}{total}{} & \psync{}{comm}{total}{}\\
    \hline
{1} & \psync{746}{0.0}{10.9}{100\%} & \psync{746}{0.0}{10.9}{100\%} & \psync{746}{0.0}{10.9}{100\%}\\
{8} & \psync{746}{34.5}{35.8}{3.81\%} & \psync{746}{1.5}{3.1}{43.81\%} & \psync{699}{1.4}{3.2}{42.09\%}\\
{16} & \psync{746}{67.0}{67.7}{1.01\%} & \psync{746}{0.3}{2.5}{26.94\%} & \psync{678}{1.5}{2.6}{25.88\%}\\
{24} & \psync{746}{98.5}{98.9}{0.46\%} & \psync{746}{0.5}{2.9}{15.96\%} & \psync{680}{1.9}{2.9}{15.52\%}\\
{32} & \psync{746}{132.6}{132.9}{0.26\%} & \psync{746}{0.6}{2.3}{14.84\%} & \psync{646}{1.7}{2.4}{14.46\%}\\
{40} & \psync{746}{169.2}{169.4}{0.16\%} & \psync{746}{0.4}{2.4}{11.38\%} & \psync{617}{1.6}{2.5}{10.93\%}\\
{48} & \psync{746}{547.5}{547.6}{0.04\%} & \psync{746}{0.6}{2.7}{8.43\%} & \psync{582}{1.9}{2.6}{8.70\%}\\
{56} & \psync{746}{650.6}{650.8}{0.03\%} & \psync{746}{0.8}{2.3}{8.48\%} & \psync{568}{1.7}{2.5}{7.68\%}\\
{64} & \psync{746}{734.6}{734.7}{0.02\%} & \psync{746}{1.1}{2.1}{8.13\%} & \psync{525}{1.8}{2.4}{7.02\%}\\
    \hlinewd{1.0pt}
  \end{tabular}
  }
}{Numerical results of \textit{luf\_cube-35937}}{}{parallel_luf_cube-35937}
Normally, the number of iterations of all synchronous algorithms should be identical independently from the number of processors for all algorithms.
The band-row versions give the same number of iterations for all numbers of processors.
In contrast, the number of iterations of the synchronous substructuring algorithm is a little different.
The effects on the number of iterations can be explained by the use of simultaneous local convergence as a stopping criterion.
In fact, the global convergence depends on the local convergence, where the local subdomain depends on the number of processors.
\includetabtbl{
  {
  \renewcommand{\arraystretch}{1.0}
  \renewcommand{\tabcolsep}{0.05cm}
  \begin{tabular}{|>\bfseries c|ccc|ccc|ccc|}
    \hlinewd{1.0pt}
    & \multicolumn{9}{c|}{\textbf{\textsc{Synchronous}}} \\
    & \multicolumn{3}{c|}{\emph{J-BANDROW}} & \multicolumn{3}{c|}{\emph{J-BANDROW-OP}} & \multicolumn{3}{c|}{\emph{J-SUBSTRUCTURING}}\\
    \hline
    \textbf{\#p (\#n)} & {\bfseries \# iter.} & {\bfseries time (s)} & {\bfseries eff.} & {\bfseries \# iter.} & {\bfseries time (s)} & {\bfseries eff.} & {\bfseries \# iter.} & {\bfseries time (s)} & {\bfseries eff.} \\
                       & \psync{}{comm}{total}{}& \psync{}{comm}{total}{} & \psync{}{comm}{total}{}\\
    \hline
{1} & \psync{2,216}{0.0}{238.4}{100\%} & \psync{2,216}{0.0}{238.4}{100\%} & \psync{2,216}{0.0}{238.4}{100\%}\\
{8} & \psync{2,216}{848.7}{882.7}{3.38\%} & \psync{2,216}{14.9}{53.9}{55.27\%} & \psync{2,132}{13.7}{48.0}{62.02\%}\\
{16} & \psync{2,216}{1,617.8}{1,634.2}{0.91\%} & \psync{2,216}{2.3}{37.1}{40.21\%} & \psync{2,098}{14.5}{32.1}{46.35\%}\\
{24} & \psync{2,216}{2,675.6}{2,686.6}{0.37\%} & \psync{2,216}{3.3}{33.8}{29.40\%} & \psync{2,067}{13.7}{27.8}{35.74\%}\\
{32} & \psync{2,216}{3,222.1}{3,229.7}{0.23\%} & \psync{2,216}{2.9}{28.7}{25.98\%} & \psync{2,041}{14.6}{24.6}{30.25\%}\\
{40} & \psync{2,216}{3,725.8}{3,731.7}{0.16\%} & \psync{2,216}{6.7}{34.5}{17.29\%} & \psync{1,995}{17.2}{27.1}{21.95\%}\\
{48} & \psync{2,216}{4,341.1}{4,345.9}{0.11\%} & \psync{2,216}{6.6}{33.6}{14.79\%} & \psync{1,980}{17.7}{25.8}{19.27\%}\\
{56} & \psync{2,216}{4,957.1}{4,961.2}{0.09\%} & \psync{2,216}{6.2}{31.2}{13.64\%} & \psync{1,907}{16.5}{24.7}{17.26\%}\\
{64} & \psync{2,216}{7,055.0}{7,058.6}{0.05\%} & \psync{2,216}{5.7}{28.0}{13.31\%} & \psync{1,889}{16.2}{23.0}{16.18\%}\\
    \hlinewd{1.0pt}
  \end{tabular}
  }
}{Numerical results of \textit{luf\_cube-274625}}{}{parallel_luf_cube-274625}

In parallel synchronous algorithms, the substructuring is always fast, in particular with a large number of processors.
The naive band-row algorithm is the slowest.
This version is slow because of a growing number of communications.
\includetabtbl{
  {
  \renewcommand{\arraystretch}{1.0}
  \renewcommand{\tabcolsep}{0.05cm}
  \begin{tabular}{|>\bfseries c|ccc|ccc|ccc|}
    \hlinewd{1.0pt}
    & \multicolumn{9}{c|}{\textbf{\textsc{Synchronous}}} \\
    & \multicolumn{3}{c|}{\emph{J-BANDROW}} & \multicolumn{3}{c|}{\emph{J-BANDROW-OP}} & \multicolumn{3}{c|}{\emph{J-SUBSTRUCTURING}}\\
    \hline
    \textbf{\#p (\#n)} & {\bfseries \# iter.} & {\bfseries time (s)} & {\bfseries eff.} & {\bfseries \# iter.} & {\bfseries time (s)} & {\bfseries eff.} & {\bfseries \# iter.} & {\bfseries time (s)} & {\bfseries eff.} \\
                       & \psync{}{comm}{total}{}& \psync{}{comm}{total}{} & \psync{}{comm}{total}{}\\
    \hline
{1} & \psync{6,234}{0.0}{159.5}{100\%} & \psync{6,234}{0.0}{159.5}{100\%} & \psync{6,234}{0.0}{159.5}{100\%}\\
{8} & \psync{6,234}{650.1}{668.2}{2.98\%} & \psync{6,234}{121.8}{142.7}{13.97\%} & \psync{6,234}{8.4}{33.8}{59.00\%}\\
{16} & \psync{6,234}{1,099.4}{1,108.5}{0.90\%} & \psync{6,234}{105.2}{133.1}{7.49\%} & \psync{6,184}{9.1}{22.9}{43.58\%}\\
{24} & \psync{6,234}{1,644.8}{1,650.9}{0.40\%} & \psync{6,234}{107.3}{137.4}{4.83\%} & \psync{6,130}{12.0}{25.7}{25.82\%}\\
{32} & \psync{6,234}{2,284.2}{2,288.7}{0.22\%} & \psync{6,234}{105.3}{136.6}{3.65\%} & \psync{6,163}{14.7}{23.7}{21.02\%}\\
{40} & \psync{6,234}{5,655.7}{5,659.4}{0.07\%} & \psync{6,234}{103.6}{152.2}{2.62\%} & \psync{5,940}{12.7}{22.8}{17.50\%}\\
{48} & \psync{6,234}{6,889.4}{6,892.5}{0.05\%} & \psync{6,234}{110.7}{151.7}{2.19\%} & \psync{6,201}{10.4}{22.9}{14.48\%}\\
{56} & \psync{6,234}{7,422.6}{7,425.3}{0.04\%} & \psync{6,234}{108.2}{150.7}{1.89\%} & \psync{6,025}{11.8}{21.7}{13.11\%}\\
{64} & \psync{6,234}{8,435.8}{8,438.5}{0.03\%} & \psync{6,234}{98.1}{149.2}{1.67\%} & \psync{6,011}{14.5}{22.6}{11.02\%}\\
    \hlinewd{1.0pt}
  \end{tabular}
  }
}{Numerical results of \textit{gravi\_hexas100x100x6\_0}}{}{parallel_gravi_hexas100x100x6_0}

\includetabtbl{
  {
  \renewcommand{\arraystretch}{1.0}
  \renewcommand{\tabcolsep}{0.035cm}
  \begin{tabular}{|>\bfseries c|ccc|ccc|ccc|}
    \hlinewd{1.0pt}
    & \multicolumn{9}{c|}{\textbf{\textsc{Synchronous}}} \\
    & \multicolumn{3}{c|}{\emph{J-BANDROW}} & \multicolumn{3}{c|}{\emph{J-BANDROW-OP}} & \multicolumn{3}{c|}{\emph{J-SUBSTRUCTURING}}\\
    \hline
    \textbf{\#p (\#n)} & {\bfseries \# iter.} & {\bfseries time (s)} & {\bfseries eff.} & {\bfseries \# iter.} & {\bfseries time (s)} & {\bfseries eff.} & {\bfseries \# iter.} & {\bfseries time (s)} & {\bfseries eff.} \\
                       & \psync{}{comm}{total}{}& \psync{}{comm}{total}{} & \psync{}{comm}{total}{}\\
    \hline
{1} & \psync{13,791}{0.0}{2,713.0}{100\%} & \psync{13,791}{0.0}{2,713.0}{100\%} & \psync{13,791}{0.0}{2,713.0}{100\%}\\
{8} & \psync{13,791}{9,508.5}{9,907.6}{3.42\%} & \psync{13,791}{5,994.3}{6,891.7}{4.92\%} & \psync{13,791}{66.0}{469.7}{72.20\%}\\
{16} & \psync{13,791}{19,564.6}{19,760.2}{0.86\%} & \psync{13,791}{6,869.0}{7,399.3}{2.29\%} & \psync{13,637}{66.0}{273.7}{61.95\%}\\
{24} & \psync{13,791}{27,329.7}{27,453.9}{0.41\%} & \psync{13,791}{7,002.2}{7,689.7}{1.47\%} & \psync{13,791}{71.9}{222.8}{50.73\%}\\
{32} & \psync{13,791}{34,939.7}{35,030.2}{0.24\%} & \psync{13,791}{7,309.4}{8,429.6}{1.01\%} & \psync{13,456}{97.7}{228.3}{37.13\%}\\
{40} & \psync{13,791}{42,060.2}{42,132.9}{0.16\%} & \psync{13,791}{7,475.9}{9,768.5}{0.69\%} & \psync{13,734}{116.4}{210.4}{32.24\%}\\
{48} & \psync{13,791}{49,218.6}{49,278.7}{0.11\%} & \psync{13,791}{6,513.5}{10,769.9}{0.52\%} & \psync{13,674}{104.9}{191.9}{29.46\%}\\
{56} & \psync{13,791}{56,170.8}{56,223.9}{0.09\%} & \psync{13,791}{6,370.3}{10,940.9}{0.44\%} & \psync{13,590}{105.5}{187.8}{25.80\%}\\
{64} & \psync{13,791}{62,913.8}{62,965.6}{0.07\%} & \psync{13,791}{7,770.8}{11,921.9}{0.36\%} & \psync{13,069}{86.6}{164.6}{25.76\%}\\
    \hlinewd{1.0pt}
  \end{tabular}
  }
}{Numerical results of \textit{gravi\_hexas100x100x6\_1}}{}{parallel_gravi_hexas100x100x6_1}

We can see in all tables that with the optimized version we considerably decrease the number of exchanges and therefore improve the execution time of the solver.
These results are even more important when the problem is large in size, as we can observe in \tabletblreft{parallel_gravi_hexas100x100x6_1}.
To overcome the lake of the synchronous algorithms, asynchronous variants can be considered, especially with a large number of processes.
However, for small problems asynchronous algorithms may be less efficient.

\section{Conclusion}
\label{sec:conclusion}
This paper evaluates and analyses the performance of parallel synchronous Jacobi algorithms by different splitting strategies including band-row splitting, band-row sparsity pattern splitting and substructuring splitting.
Matrices arising from the finite element discretization of diverse equations (3D Laplace equation, 3D gravitational potential) with different meshes are used as a testbed.
The results clearly show the interest of the substructuring methods compare to band-row parallel synchronous algorithms to solve sparse large linear systems.

\bibliography{bib/zotero_20160504}
\bibliographystyle{abbrv}

\end{document}